\documentclass[final,1p,times]{elsarticle}

\usepackage{amssymb}
 \usepackage{amsthm}
\usepackage{amscd}
\usepackage{amsmath}
\usepackage{amsfonts}
\usepackage{amssymb}
\usepackage{graphicx}
\newtheorem{theorem}{Theorem}

\newtheorem{lemma}[theorem]{Lemma}

\newtheorem{remark}[theorem]{Remark}
\usepackage{mathrsfs}
\usepackage{titletoc}

\newcommand{\ra}{\rightarrow}
\newcommand{\p}{\partial}
\newcommand{\f}{\frac}

\newcommand{\be}{\begin{equation}}
\renewcommand{\ra}{\rightarrow}
\newcommand{\ee}{\end{equation}}
\newcommand{\bea}{\begin{eqnarray}}
\newcommand{\eea}{\end{eqnarray}}
\newcommand{\bna}{\begin{eqnarray*}}
\newcommand{\ena}{\end{eqnarray*}}

\renewcommand{\le}{\left}
\newcommand{\ri}{\right}

\journal{***}

\begin{document}

\begin{frontmatter}
\title{A gradient estimate for positive functions on graphs}

\author{Yong Lin}
 \ead{linyong01@ruc.edu.cn}

\author{Shuang Liu}
 \ead{cherrybu@ruc.edu.cn}

\author{Yunyan Yang}
 \ead{yunyanyang@ruc.edu.cn}
\address{ Department of Mathematics,
Renmin University of China, Beijing 100872, P. R. China}

\begin{abstract}

  We derive a gradient estimate for positive functions, in particular for positive solutions to
  the heat equation, on finite or locally finite graphs. Unlike the well known Li-Yau estimate,
  which is based on the maximum principle, our estimate follows from the graph structure of
  the gradient form and the Laplacian operator. Though our assumption on graphs is slightly stronger than that of
  Bauer, Horn, Lin, Lippner, Mangoubi, and Yau (J. Differential Geom. 99 (2015) 359-405), our estimate can be easily applied
  to nonlinear differential equations, as well as differential inequalities.
  As applications, we estimate the greatest lower bound of Cheng's eigenvalue and an upper bound of
  the minimal heat kernel, which is recently studied by Bauer, Hua and Yau (Preprint, 2015)
  by the Li-Yau estimate. Moreover, generalizing an earlier result of Lin and Yau (Math. Res. Lett. 17 (2010) 343-356),
  we derive a lower bound of nonzero eigenvalues by our gradient estimate.
\end{abstract}


\end{frontmatter}

\titlecontents{section}[0mm]
                       {\vspace{.2\baselineskip}}
                       {\thecontentslabel~\hspace{.5em}}
                        {}
                        {\dotfill\contentspage[{\makebox[0pt][r]{\thecontentspage}}]}
\titlecontents{subsection}[3mm]
                       {\vspace{.2\baselineskip}}
                       {\thecontentslabel~\hspace{.5em}}
                        {}
                       {\dotfill\contentspage[{\makebox[0pt][r]{\thecontentspage}}]}

\setcounter{tocdepth}{2}


\section{Introduction}

Let $G=(V,E)$ be a finite or locally finite graph, where $V$ denotes the vertex set and $E$ denotes the edge set.
For any edge $xy\in E$, we assume its weight $w_{xy}>0$. The degree of $x\in V$ is defined
as ${\rm deg}(x)=\sum_{y\sim x}w_{xy}$, here and throughout this paper we write $y\sim x$ if
$xy\in E$. Let $\mu:V\ra \mathbb{R}$ be a finite measure. Then the $\mu$-Laplacian (or Laplacian for short) on
$G$ is defined as
$$\Delta f(x)=\f{1}{\mu(x)}\sum_{y\sim x}w_{xy}(f(y)-f(x)).$$
The associated gradient form reads
$$2\Gamma(f,g)(x)=\f{1}{\mu(x)}\sum_{y\sim x}w_{xy}(f(y)-f(x))(g(y)-g(x)).$$
Write $\Gamma(f)=\Gamma(f,f)$. Denote
\be\label{cond1}D_\mu=\sup_{x\in V}\f{{\deg}(x)}{\mu(x)},\quad d=\sup_{x\in V,\,xy\in E}\f{\mu(x)}{w_{xy}}.\ee
The main result in this paper is the following gradient estimate.

\begin{theorem}\label{Theorem 1}
 Let $G=(V,E)$ be a finite or locally finite graph.  Suppose that
 \be\label{cd1}D_\mu<+\infty,\quad d<+\infty,\ee
 where $D_\mu$ and $d$ are defined as in (\ref{cond1}).
 Then
   for any positive function $u: V\ra \mathbb{R}$, there holds
   \be\label{gradient1}\f{\sqrt{2\Gamma(u)}}{u}\leq \sqrt{d}\f{\Delta u}{u}+\sqrt{d}D_\mu+\sqrt{D_\mu}.\ee
   Several special cases are listed below:\\

   \noindent $(i)$ If $u$ is a positive solution to the differential inequality $\Delta u-qu\leq 0$ on $V$, where $q:V\ra \mathbb{R}$
is a function, then there holds
$$\f{\sqrt{2\Gamma(u)}}{u}-\sqrt{d}q\leq\sqrt{d}D_\mu+\sqrt{D_\mu}.$$

\noindent $(ii)$ If $u$ is a positive solution to the differential inequality  $\Delta u-h u^\alpha\leq 0$, where $\alpha\in \mathbb{R}$,
and $h:V\ra \mathbb{R}$ is a function, then there holds
$$\f{\sqrt{2\Gamma(u)}}{u}-\sqrt{d}hu^{\alpha-1}\leq\sqrt{d}D_\mu+\sqrt{D_\mu}.$$

\noindent $(iii)$ If $u$ is a positive solution to the differential inequality  $\Delta u-\p_tu\leq qu$, where $q:V\times\mathbb{R}
\ra \mathbb{R}$
is a function, then there holds
$$\f{\sqrt{2\Gamma(u)}}{u}-\sqrt{d}\f{\p_tu}{u}-\sqrt{d}q\leq \sqrt{d}D_\mu+\sqrt{D_\mu}.$$

\noindent $(iv)$ If $u$ is a positive solution to the differential inequality
$\Delta u-\p_tu+au\log u\leq 0$, where $a\in\mathbb{R}$ is a constant, then there holds
$$\f{\sqrt{2\Gamma(u)}}{u}-\sqrt{d}\f{\p_tu}{u}-\sqrt{d}a\log u\leq \sqrt{d}D_\mu+\sqrt{D_\mu}.$$
\end{theorem}

\begin{remark}
For the corresponding partial differential equations on complete Riemannian manifolds, $(i)-(iv)$
were extensively studied, see for examples \cite{LiYau1,LiYau2,LiJ,Negrin,Ma,Yang-1,Yang-2} and the references there in.
\end{remark}


At least two points can be seen from Theorem \ref{Theorem 1}: One is that (\ref{gradient1}) is a global estimate;
The other is that (\ref{gradient1}) can be easily applied to nonlinear elliptic or parabolic equations, as well as
differential inequalities. We now analyze the assumption (\ref{cd1}), whch is equivalent to
\be\label{cd2}\sup_{x\in V}\sharp\le\{y|y\sim x\ri\}<+\infty,\quad
0<\inf_{x\in V,\,y\sim x}\f{\mu(x)}{w_{xy}}\leq \sup_{x\in V,\,y\sim x}\f{\mu(x)}{w_{xy}}<+\infty,\ee
where $\sharp\le\{y|y\sim x\ri\}$ stands for the number of $y\in V$ which is adjacent to $x$. In fact,
suppose (\ref{cd1}) holds. Then $\sharp\le\{y|y\sim x\ri\}\leq D_\mu d$ and $\f{1}{D_\mu}\leq \f{\mu(x)}{w_{xy}}\leq d$ for
any $y\sim x$. Hence (\ref{cd2}) holds. Conversely, if (\ref{cd2}) holds, we have
$$\f{{\rm deg}(x)}{\mu(x)}=\f{\sum_{y\sim x}w_{xy}}{\mu(x)}\leq\f{\sharp\{y|y\sim x\}}
{\inf_{x\in V,\,y\sim x}\f{\mu(x)}{w_{xy}}}.$$
Then (\ref{cd1}) follows immediately. If we replace (\ref{cd1}) by
(\ref{cd2}) in Theorem \ref{Theorem 1}, then the gradient estimate (\ref{gradient1})
would be
\be\label{gradient1'}\f{\sqrt{2\Gamma(u)}}{u}\leq \sqrt{b}\f{\Delta u}{u}+\sqrt{b}\le(\f{N}{a}+\sqrt{\f{N}{a b}}\ri),\ee
where $N=\sup_{x\in V}\sharp\le\{y|y\sim x\ri\}$, $a=\inf_{x\in V,\,y\sim x}\f{\mu(x)}{w_{xy}}$
and $b=\sup_{x\in V,\,y\sim x}\f{\mu(x)}{w_{xy}}$.
Note that the essential assumption of \cite{BHLLY}  is
\be\label{BLLY}D_\mu<+\infty,\quad D_w=\sup_{x\in V,\,y\sim x}\f{{\rm deg}(x)}{w_{xy}}<+\infty.\ee
It is easy to see that (\ref{BLLY}) is slightly weaker than (\ref{cd1}). All gradient estimates in \cite{BHLLY} are about
$\sqrt{u}$,  where $u$ is a positive solution to a parabolic equation on $V$. Note that for any positive function $u$
\bea\nonumber
2\Gamma(\sqrt{u})(x)&=&\f{1}{\mu(x)}\sum_{y\sim x}w_{xy}\le(\sqrt{u}(y)-\sqrt{u}(x)\ri)^2\\\nonumber
&\leq&\le(\sum_{y\sim x}\f{w_{xy}}{\mu(x)}\ri)^{1/2}\le(\f{1}{\mu(x)}\sum_{y\sim x}w_{xy}
\le(\sqrt{u}(y)-\sqrt{u}(x)\ri)^4\ri)^{1/2}\\\nonumber
&\leq&\le(\f{{\rm deg}(x)}{\mu(x)}\ri)^{1/2}\le(\f{1}{\mu(x)}\sum_{y\sim x}w_{xy}
\le({u}(y)-{u}(x)\ri)^2\ri)^{1/2}\\\label{r-l}
&\leq&\sqrt{D_\mu}\sqrt{2\Gamma(u)(x)}.
\eea
If $u:V\times [0,+\infty)\ra \mathbb{R}$ is a positive solution to the parabolic equation $\Delta u-\p_tu-qu=0$
on $V\times[0,+\infty)$, we conclude an analog of (\cite{BHLLY}, Theorem 4.10) by combining (\ref{r-l}) with Theorem \ref{Theorem 1},
$$\f{\Gamma(\sqrt{u})}{u}-\sqrt{D_\mu d}\f{\p_t\sqrt{u}}{\sqrt{u}}-\sqrt{D_\mu d}\f{q}{2}
\leq \f{D_\mu(\sqrt{D_\mu d}+1)}{2}.$$

As an application of Case $(iii)$ of Theorem \ref{Theorem 1}, we state the following Harnack inequality.

\begin{theorem}{\label{Theorem2}} Let $G=(V,E)$ be a finite or locally finite graph satisfying (\ref{cd1}).
Moreover $\mu_{\max}=\sup_{x\in V}\mu(x)<+\infty$ and $w_{\min}=\inf_{x\in V,\,y\sim x} w_{xy}>0$.
Assume $u:V\times(-\infty,+\infty)\ra \mathbb{R}$ is a positive solution to the heat inequality
$\Delta u-\p_tu\leq qu$, where $q: V\times (-\infty,+\infty)\ra \mathbb{R}$ is a function.
Then for any $(x,T_1)$ and $(y,T_2)$, $T_1<T_2$, we have
\bna
u(x,T_1)&\leq& u(y,T_2)\exp\le\{\le(D_\mu+\sqrt{\f{D_\mu}{d}}\ri)(T_2-T_1)+\f{\le({\rm dist}(x,y)\ri)^2}{T_2-T_1}
\sqrt{\f{d\mu_{\max}}{w_{\min}}}\ri.\\
&&\le.+\min \sum_{k=0}^{\ell-1}\le(\int_{t_k}^{t_{k+1}}q(x_k,t)dt+\f{\ell^2}{(T_2-T_1)^2}\int_{t_k}^{t_{k+1}}(t-t_k)^2
(q(x_{k+1},t)-q(x_k,t))dt\ri)\ri\},
\ena
where the minimum takes over all shortest paths $x=x_0, x_1,\cdots, x_\ell=y$ connecting $x$ and $y$,
and $t_k=T_1+k(T_2-T_1)/\ell$, $k=0,1,\cdots,\ell$. In particular, if there exists some constant $C_0$ such that
$|q(x,t)|\leq C_0$ for all $(x,t)$, then for any $x,y\in V$ and $T_1<T_2$,
there holds
\be\label{Har}u(x,T_1)\leq u(y,T_2)\exp\le\{\le(D_\mu+\sqrt{\f{D_\mu}{d}}+\f{5}{3}C_0\ri)(T_2-T_1)+\f{\le({\rm dist}(x,y)\ri)^2}{T_2-T_1}
\sqrt{\f{d\mu_{\max}}{w_{\min}}}\ri\}.\ee
\end{theorem}

In a recent work of F. Bauer, B. Hua and S. T. Yau \cite{BHY}, the Li-Yau inequality on graphs, which is due to
F. Bauer, P. Horn, Y. Lin, G. Lipper, D. Mangoubi and S. T. Yau \cite{BHLLY}, is applied to Liouville type theorems and eigenvalue
estimates. Moreover a DGG lemma \cite{Gaf59,Dav92,Gri94} concerning the minimal heat kernel is established on graphs, and it is used together with the Li-Yau
inequality to estimate the upper bound of the minimal heat kernel. Our gradient estimate can be used instead of the Li-Yau estimate
in \cite{BHY}. Using Theorem \ref{Theorem 1}, we can estimate the greatest lower bound of the $\ell^2$-spectrum known
as Cheng's eigenvalue \cite{Che75}.

\begin{theorem}\label{Theorem 1.2}
Let $G=(V,E)$ be a locally finite graph satisfying (\ref{cd1}) and $\lambda^\ast$ be the greatest lower bound of the $\ell^2$-spectrum
of the graph Laplacian $\Delta$. Then we have $\lambda^\ast\leq D_\mu+\sqrt{D_\mu}/\sqrt{d}$.
\end{theorem}
While Theorem \ref{Theorem2} can be used to get an analog of (\cite{BHY}, Theorem 1.2).

\begin{theorem}\label{Theorem 1.3}
Let $G=(V,E)$ be a finite or locally finite graph satisfying (\ref{cd1}). Moreover $\mu_{\max}<+\infty$ and
$w_{\min}>0$. Let $\lambda^\ast$ be the greatest lower bound of the $\ell^2$-spectrum
of the graph Laplacian $\Delta$.
Given any $\epsilon>0$, $0<\gamma\leq1$, $\beta>0$. Let $P_t(x,y)$ be the minimum heat kernel of $G$. Then there exist
positive constants
$C_1(\beta,\gamma,D_\mu)$ and $C_2(\epsilon,\beta,\gamma,D_\mu,d,\mu_{\max},w_{\min})$
such that for any $x,y\in V$ and
$t\geq \max\{\beta d(x,y),1\}$,
$$P_t(x,y)\leq \f{\exp\le(-(1-\gamma)\lambda^\ast t\ri)}{\sqrt{{\rm Vol}(B_x(\sqrt{t})){\rm Vol}(B_y(\sqrt{t}))}}
\exp\le\{C_2\sqrt{t}-C_1\f{({\rm dist}(x,y))^2}{4(1+2\epsilon)t}\ri\}.$$
\end{theorem}

Finally we remark that Theorem 1 can also be used to estimate a lower bound of nonzero eigenvalues of
the Laplacian on
finite connected graphs. Precisely we have an analog of (\cite{LinYau}, Theorem 1.8), namely
\begin{theorem}\label{lowerbd}
Let $G=(V,E)$ be a finite connected graph, $D_\mu$ and $d$ be defined as in (\ref{cond1}), and
$D$ be its diameter. Moreover we assume $w_{x,y}=w_{yx}$ for all $y\sim x$ and all $x\in V$.
Suppose that $\lambda$ is a nonzero eigenvalue of $-\Delta$. Then there holds
\be\label{lo-b}\lambda\geq\f{1}{Dd\le(\exp\le\{1+Dd\le(D_\mu+\sqrt{\f{D_\mu}{d}}\ri)\ri\}-1\ri)}.\ee
\end{theorem}

\begin{remark}
If $\mu(x)=\deg(x)=\sum_{y\sim x}w_{xy}$, we have $D_\mu=1$, and whence (\ref{lo-b}) becomes
$$\lambda\geq\f{1}{Dd\le(\exp\le\{1+Dd\le(1+\sqrt{\f{1}{d}}\ri)\ri\}-1\ri)}.$$
 We refer the reader to \cite{BCG,FanChung} for earlier estimates in terms of the volume of the graph $G$.
\end{remark}

 Let us describe the method. The proof of Theorem \ref{Theorem 1} is based on the positivity of the average of
 $u$, i.e. $\f{1}{\mu(x)}\sum_{y\sim x}w_{xy}u(y)$, and
 its relation with $\Delta u$. To prove Theorem \ref{Theorem2}, we follow \cite{BHLLY} and thereby closely
 follow \cite{LiYau2}. While the proof of Theorems \ref{Theorem 1.2}, \ref{Theorem 1.3} and \ref{lowerbd}
 is adapted from \cite{BHY} and \cite{LinYau} respectively.

The remaining part of this paper is organized as follows. In Section 2, we prove the gradient estimate, Theorem \ref{Theorem 1}.
In Section 3, we prove the corresponding Harnack inequality,
Theorem \ref{Theorem2}. Finally Theorems \ref{Theorem 1.2}, \ref{Theorem 1.3} and \ref{lowerbd} are proved in Section 4.

\section{Gradient estimate}
In this section, we prove Theorem \ref{Theorem 1} by using a very simple method.\\

{\it Proof of Theorem \ref{Theorem 1}}.
 Special cases $(i)-(iv)$ are immediate consequences of (\ref{gradient1}).
 Hence it suffices to prove (\ref{gradient1}).
Since $u>0$, we have by definition of $\Gamma(u)$,
\bna
2\Gamma(u)(x)&=&\f{1}{\mu(x)}\sum_{y\sim x}w_{xy}(u(y)-u(x))^2\\
&\leq&\f{1}{\mu(x)}\sum_{y\sim x}w_{xy}u^2(y)+\f{1}{\mu(x)}\sum_{y\sim x}w_{xy}u^2(x)\\
&=&\sum_{y\sim x}\f{\mu(x)}{w_{xy}}\le(\f{w_{xy}}{\mu(x)}u(y)\ri)^2+u^2(x)\f{{\deg}(x)}{\mu(x)}\\
&\leq&d\le(\sum_{y\sim x}\f{w_{xy}}{\mu(x)}u(y)\ri)^2+D_\mu u^2(x).
\ena
Noting that
$$\sum_{y\sim x}\f{w_{xy}}{\mu(x)}u(y)=\Delta u(x)+u(x)\f{{\deg}(x)}{\mu(x)}$$
and using an elementary inequality $\sqrt{a^2+b^2}\leq a+b$, $\forall a, b\geq 0$,
we get
\bna
\sqrt{2\Gamma(u)(x)}&\leq&\sqrt{d}\sum_{y\sim x}\f{w_{xy}}{\mu(x)}u(y)+\sqrt{D_\mu}u(x)\\
&=&\sqrt{d}\le(\Delta u(x)+u(x)\f{{\deg}(x)}{\mu(x)}\ri)+\sqrt{D_\mu}u(x)\\
&\leq&\sqrt{d}\Delta u(x)+(\sqrt{D_\mu}+\sqrt{d}D_\mu)u(x).
\ena
This leads to (\ref{gradient1}) and thus ends the proof of the theorem. $\hfill\Box$\\

\section{Harnack inequality}

In this section, following the lines of \cite{BHLLY,LiYau2}, we prove a Harnack inequality for
positive solution to the parabolic inequality $\Delta u-\p_tu\leq qu$ by using $(iii)$ of Theorem \ref{Theorem 1}.\\

{\it Proof of Theorem \ref{Theorem2}}. Let $u$ be a positive solution to the inequality $\Delta u-\p_tu\leq qu$.
 By $(iii)$ of Theorem \ref{Theorem 1}, we have
\be\label{dtlog}-\p_t\log u\leq D_\mu+\f{\sqrt{D_\mu}}{\sqrt{d}}+q-\f{1}{\sqrt{d}}\f{\sqrt{2\Gamma(u)}}{u}.\ee
We distinguish two cases to proceed.\\

{\bf Case 1}. $x\sim y$.

For any $s\in [T_1,T_2]$, we have by (\ref{dtlog}) that
\bea\nonumber
\log u(x,T_1)-\log u(y,T_2)&=&\log\f{u(x,T_1)}{u(x,s)}+\log\f{u(x,s)}{u(y,s)}+
\log\f{u(y,s)}{u(y,T_2)}\\\nonumber
&=&-\int_{T_1}^s\p_t\log u(x,t)dt+\log\f{u(x,s)}{u(y,s)}-
\int_s^{T_2}\p_t\log u(y,t)dt\\\nonumber
&\leq&\le(D_\mu+\f{\sqrt{D_\mu}}{\sqrt{d}}\ri)(T_2-T_1)+\int_{T_1}^sq(x,t)dt+
\int_s^{T_2}q(y,t)dt
\\&&-\f{1}{\sqrt{d}}\le(\int_{T_1}^s\f{\sqrt{2\Gamma(u)(x,t)}}{u(x,t)}dt+
\int_s^{T_2}\f{\sqrt{2\Gamma(u)(y,t)}}{u(y,t)}dt\ri)\nonumber\\&&+\log\f{u(x,s)}{u(y,s)}.
\label{temp}\eea
We estimate the above terms respectively. Obviously
\be\label{pos}\int_{T_1}^{s}\f{\sqrt{2\Gamma(u)(x,t)}}{u(x,t)}dt\geq 0.\ee
Since
\bna
2\Gamma(u)(y,t)&=&\f{1}{\mu(y)}\sum_{z\thicksim y}w_{yz}(u(z,t)-u(y,t))^2\\
&\geq&\f{w_{\min}}{\mu_{\max}}(u(x,t)-u(y,t))^2,
\ena
we get
\be\label{2gam}
-\f{1}{\sqrt{d}}
\int_s^{T_2}\f{\sqrt{2\Gamma(u)(y,t)}}{u(y,t)}dt
\leq -\sqrt{\f{w_{\min}}{d\mu_{\max}}}\int_s^{T_2}\le|\f{u(x,t)}{u(y,t)}-1\ri|dt.
\ee
Using an elementary inequality $\log r\leq \sqrt{|r-1|}, \forall r>0$, we have
\be\label{log}\log\f{u(x,s)}{u(y,s)}\leq \psi(x,y,s),\ee
where
$$\psi(x,y,s)=\sqrt{\le|\f{u(x,s)}{u(y,s)}-1\ri|}.$$
Inserting (\ref{pos}), (\ref{2gam}), (\ref{log}) into (\ref{temp}), and using
(\cite{BHLLY}, Lemma 5.3), we obtain
\bna
\log u(x,T_1)-\log u(y,T_2)&\leq&\le(D_\mu+\sqrt{\f{D_\mu}{d}}\ri)(T_2-T_1)+\psi(x,y,s)\\
&&-\sqrt{\f{w_{\min}}{d\mu_{\max}}}\int_s^{T_2}\psi^2(x,y,t)dt
+\int_{T_1}^sq(x,t)dt+\int_{s}^{T_2}q(y,t)dt\\
&\leq&\le(D_\mu+\sqrt{\f{D_\mu}{d}}\ri)(T_2-T_1)+\f{1}{T_2-T_1}\sqrt{\f{d\mu_{\max}}{w_{\min}}}\\
&&+\int_{T_1}^{T_2}q(x,t)dt+\f{1}{(T_2-T_1)^2}\int_{T_1}^{T_2}(t-T_1)^2(q(y,t)-q(x,t))dt.
\ena

{\bf Case 2}. $x$ is not adjacent to $y$.\\

Assume ${\rm dist}(x,y)=\ell$. Take a shortest path $x=x_0, x_1,\cdots,x_{\ell}=y$. Let $T_1=t_0<t_1<\cdots
<t_\ell=T_2$, $t_k=t_{k-1}+(T_2-T_1)/\ell$, $k=1,\cdots,\ell$. By the result of Case 1, we have
\bna
&&\log u(x,T_1)-\log u(y,T_2)=\sum_{k=0}^{\ell-1}\le(\log u(x_k,t_k)-\log u(x_{k+1},t_{k+1})\ri)\\
&&\leq\sum_{k=0}^{\ell-1}\le(\le(D_\mu+\sqrt{\f{D_\mu}{d}}\ri)(t_{k+1}-t_k)+\f{1}{t_{k+1}-t_k}\sqrt{\f{d\mu_{\max}}{w_{\min}}}\ri)\\
&&\quad+\sum_{k=0}^{\ell-1}\le(\int_{t_k}^{t_{k+1}}q(x_k,t)dt+\f{1}{(t_{k+1}-t_k)^2}\int_{t_k}^{t_{k+1}}(t-t_k)^2
(q(x_{k+1},t)-q(x_k,t))dt\ri)\\
&&\leq\le(D_\mu+\sqrt{\f{D_\mu}{d}}\ri)(T_2-T_1)+\f{\ell^2}{T_2-T_1}\sqrt{\f{d\mu_{\max}}{w_{\min}}}\\
&&\quad+\sum_{k=0}^{\ell-1}\le(\int_{t_k}^{t_{k+1}}q(x_k,t)dt+\f{\ell^2}{(T_2-T_1)^2}\int_{t_k}^{t_{k+1}}(t-t_k)^2
(q(x_{k+1},t)-q(x_k,t))dt\ri).
\ena
Therefore we conclude
\bna\log u(x,T_1)-\log u(y,T_2)&\leq& \le(D_\mu+\sqrt{\f{D_\mu}{d}}\ri)(T_2-T_1)+\f{\le({\rm dist}(x,y)\ri)^2}{(T_2-T_1)}
\sqrt{\f{d\mu_{\max}}{w_{\min}}}\\&&\quad +\min \mathscr{F}(q)(x,y,T_1,T_2), \ena
where
$$\mathscr{F}(q)(x,y,T_1,T_2)=\sum_{k=0}^{\ell-1}\le(\int_{t_k}^{t_{k+1}}q(x_k,t)dt+\f{\ell^2}{(T_2-T_1)^2}\int_{t_k}^{t_{k+1}}(t-t_k)^2
(q(x_{k+1},t)-q(x_k,t))dt\ri)$$
and the minimum takes over all shortest paths connecting $x$ and $y$. Hence the first assertion of the theorem follows
immediately.

Moreover, if $|q(x,t)|\leq C_0$ for all $(x,t)$, then we have
$$\sum_{k=0}^{\ell-1}\int_{t_k}^{t_{k+1}}q(x_k,t)dt\leq C_0(T_2-T_1)$$
and
$$\sum_{k=0}^{\ell-1}\f{\ell^2}{(T_2-T_1)^2}\int_{t_k}^{t_{k+1}}(t-t_k)^2
(q(x_{k+1},t)-q(x_k,t))dt\leq \f{2C_0}{3}(T_2-T_1).$$
This gives the desired result and the proof of the theorem is completed. $\hfill\Box$

\section{Further applications of the gradient estimate}

In this section, as applications of Theorem \ref{Theorem 1}, we prove Theorems \ref{Theorem 1.2}, \ref{Theorem 1.3} and
\ref{lowerbd}. For the proof of Theorems \ref{Theorem 1.2} and \ref{Theorem 1.3}, we follow the lines of \cite{BHY}, the
essential difference
is that we use Theorem \ref{Theorem 1} instead of the Li-Yau estimate \cite{BHLLY}. While the proof of Theorem \ref{lowerbd}
is an adaptation of \cite{LinYau}. For reader's convenience, we give the details here.\\

{\it Proof of Theorem \ref{Theorem 1.2}}. Let $\lambda^\ast$ be the greatest lower bound of Cheng's eigenvalues.
By a result of S. Haeseler and M. Keller (\cite{HK11}, Theorem 3.1),
if $\lambda<\lambda^\ast$, then there would be a positive solution $u$ to $\Delta u=-\lambda u$. We conclude from
Case $(i)$ of Theorem \ref{Theorem 1} that
$$\f{\sqrt{2\Gamma(u)}}{u}+\sqrt{d}\lambda\leq \sqrt{d}D_\mu+\sqrt{D_\mu}.$$
Hence $\lambda\leq D_\mu+\sqrt{D_\mu/d}$. Since $\lambda$ is arbitrary, we obtain $\lambda^\ast\leq D_\mu+\sqrt{D_\mu/d}$.
$\hfill\Box$\\

To prove Theorem \ref{Theorem 1.3}, we need the following DGG lemma on graphs (\cite{BHY}, Theorem 1.1).

\begin{lemma}\label{DGGlemma}
Let $P_t(x,y)$ be the minimal heat kernel of the graph $G=(V,E)$. Then for any $\beta>0$ and $0<\gamma\leq 1$,
there exists a constant $C_1$ depending only on $\beta,\gamma$ and $D_\mu$ such that
for any subsets $B_1, B_2\subset G$, $t\geq \max\{\beta {\rm dist}(B_1,B_2),1\}$,
$$\sum_{x\in B_1}\sum_{y\in B_2}P_t(x,y)\mu(x)\mu(y)\leq e^{-(1-\gamma)\lambda^\ast t}
\sqrt{{\rm Vol}(B_1){\rm Vol}(B_2)}\exp\le(-C_1\f{\le({\rm dist}(B_1, B_2)\ri)^2}{4t}\ri).$$
\end{lemma}

{\it Proof of Theorem \ref{Theorem 1.3}}. Fix $x, y\in V$, $\delta>0$, $T_1=t$ and $T_2=(1+\delta)t$.
Applying the Harnack inequality, Theorem \ref{Theorem2},
to the minimal heat kernel $P_t(x,y)$,
\bna
P_t(x,y)&\leq&P_{(1+\delta)t}(x^\prime,y)\exp\le\{\le(D_\mu+\sqrt{\f{D_\mu}{d}}\ri)\delta t+
\f{({\rm dist}(x,x^\prime))^2}{\delta t}\sqrt{\f{d\mu_{\max}}{w_{\min}}}\ri\}\\
&\leq&P_{(1+\delta)t}(x^\prime,y)\exp\le\{\le(D_\mu+\sqrt{\f{D_\mu}{d}}\ri)\delta t+
\f{1}{\delta }\sqrt{\f{d\mu_{\max}}{w_{\min}}}\ri\},\quad\forall x^\prime\in B_x(\sqrt{t}).
\ena
Integrating the above inequality on $B_x(\sqrt{t})$ with respect to $x^\prime$, we have
\be\label{vx}
{\rm Vol}(B_x(\sqrt{t}))P_t(x,y)\leq\exp\le\{\le(D_\mu+\sqrt{\f{D_\mu}{d}}\ri)\delta t+
\f{1}{\delta }\sqrt{\f{d\mu_{\max}}{w_{\min}}}\ri\}
 \sum_{x^\prime\in B_x(\sqrt{t})}\mu(x^\prime)P_{(1+\delta)t}(x^\prime,y).
\ee
Note that
$h(y,s)=\sum_{x^\prime\in B_x(\sqrt{t})}\mu(x^\prime)P_{s}(x^\prime,y)$ is also a positive solution to the heat equation.
Applying again the Harnack inequality, Theorem \ref{Theorem2}, to $h(y,s)$ with $T_1=(1+\delta)t$ and $T_2=(1+2\delta)t$,
we have
$${\rm Vol}(B_y(\sqrt{t}))h(y,(1+\delta)t)\leq \exp\le\{\le(D_\mu+\sqrt{\f{D_\mu}{d}}\ri)\delta t+
\f{1}{\delta }\sqrt{\f{d\mu_{\max}}{w_{\min}}}\ri\}
 \sum_{y^\prime\in B_y(\sqrt{t})}\mu(y^\prime)h(y^\prime,(1+2\delta)t).$$
This together with (\ref{vx}) implies that
\bna
P_t(x,y)&\leq&\exp\le\{2\le(D_\mu+\sqrt{\f{D_\mu}{d}}\ri)\delta t+
\f{2}{\delta }\sqrt{\f{d\mu_{\max}}{w_{\min}}}\ri\}\\
&&\f{1}{{\rm Vol}(B_x(\sqrt{t})){\rm Vol}(B_y(\sqrt{t}))}\sum_{x^\prime\in B_x(\sqrt{t})}\sum_{y^\prime\in B_y(\sqrt{t})}
\mu(x^\prime)\mu(y^\prime)P_{(1+2\delta)t}(x^\prime,y^\prime).
\ena

Let $t\geq \max\{\beta {\rm dist}(x,y),1\}$. Obviously $t\geq \f{1}{1+2\delta}\max
\{\beta{\rm dist}(B_x(\sqrt{t}),B_y(\sqrt{t})),1\}$. Let $\gamma$, $0<\gamma\leq 1$, be fixed.
It follows from Lemma \ref{DGGlemma} that
there exists a constant $C_1$ depending only on $\gamma$, $\beta$ and $D_\mu$ such that
\bea\nonumber
P_t(x,y)&\leq& \exp\le\{2\le(D_\mu+\sqrt{\f{D_\mu}{d}}\ri)\delta t+
\f{2}{\delta }\sqrt{\f{d\mu_{\max}}{w_{\min}}}\ri\}\f{1}
{\sqrt{{\rm Vol}(B_x(\sqrt{t})){\rm Vol}(B_y(\sqrt{t}))}}\\\label{pp}
&&\exp\le\{-(1-\gamma)\lambda^\ast(1+2\delta)t-C_1
\f{\le({\rm dist}(B_x(\sqrt{t}),B_y(\sqrt{t}))\ri)^2}{4(1+2\delta)t}\ri\}.
\eea
If ${\rm dist}(x,y)>2\sqrt{t}$, then ${\rm dist}(B_x(\sqrt{t}),B_y(\sqrt{t}))\geq
{\rm dist}(x,y)-2\sqrt{t}$ and thus
\be\label{gg}
\f{\le({\rm dist}(B_x(\sqrt{t}),B_y(\sqrt{t})\ri)^2}{4(1+2\delta)t}\geq\f{({\rm dist}(x,y))^2}{4(1+4\delta)t}
-\f{1}{2\delta}.
\ee
It is easy to see that (\ref{gg}) still holds if ${\rm dist}(x,y)\leq2\sqrt{t}$. Inserting (\ref{gg}) into (\ref{pp}),
we have
\bea\nonumber
P_t(x,y)&\leq&\f{1}
{\sqrt{{\rm Vol}(B_x(\sqrt{t})){\rm Vol}(B_y(\sqrt{t}))}}\exp\le\{2\le(D_\mu+\sqrt{\f{D_\mu}{d}}\ri)\delta t+
\f{2}{\delta }\sqrt{\f{d\mu_{\max}}{w_{\min}}}+\f{C_1}{2\delta}\ri\}\\\label{fff}
&&\exp\le\{-(1-\gamma)\lambda^\ast(1+2\delta)t-C_1
\f{\le({\rm dist}(x,y)\ri)^2}{4(1+4\delta)t}\ri\}.
\eea
Note that $t\geq 1$. Choosing $2\delta=\epsilon/\sqrt{t}$ in (\ref{fff}), we obtain for  $t\geq \max\{\beta {\rm dist}(x,y),1\}$,
\bna
P_t(x,y)&\leq&\f{1}
{\sqrt{{\rm Vol}(B_x(\sqrt{t})){\rm Vol}(B_y(\sqrt{t}))}}\exp\le\{\sqrt{t}\le(D_\mu\epsilon+\sqrt{\f{D_\mu}{d}}\epsilon+
\f{4}{\epsilon}\sqrt{\f{d\mu_{\max}}{w_{\min}}}+\f{C_1}{\epsilon}\ri)\ri\}\\\label{fi}
&&\exp\le\{-(1-\gamma)\lambda^\ast t-C_1
\f{\le({\rm dist}(x,y)\ri)^2}{4(1+2\epsilon)t}\ri\}.
\ena
Denoting $C_2=D_\mu\epsilon+\sqrt{\f{D_\mu}{d}}\epsilon+
\f{4}{\epsilon}\sqrt{\f{d\mu_{\max}}{w_{\min}}}+\f{C_1}{\epsilon}$, we finish the proof of the theorem.
$\hfill\Box$\\

{\it Proof of Theorem \ref{lowerbd}}. Note that $\int_V\Delta ud\mu=0$ and that if
$-\Delta u=\lambda u$, then $-\Delta (cu)=\lambda cu$ for any constant $c\in\mathbb{R}$.
  We can assume $-\Delta u=\lambda u$ with $\sup u=1$ and $\inf u<0$.
  Take $x_1, x_\ell\in G$ such that $u(x_1)=\sup u=1$, $u(x_n)=\inf u<0$,
$x_1x_2\cdots x_\ell$ be the shortest path connecting $x_1$ and $x_\ell$, where $(x_i,x_{i+1})\in E$. Then
$\ell\leq D$. For any $\beta>1$, note that
\bea\nonumber
\f{|u(x_i)-u(x_{i+1})|}{\beta-u(x_i)}&\leq&
\f{\sqrt{\f{1}{\mu(x_i)}\sum_{x_iy\in E}\f{\mu(x_i)}{w_{x_iy}}w_{x_iy}(u(x_i)-u(y))^2}}
{\beta-u(x_i)}\\\label{uu}
&\leq&\sqrt{d}\f{\sqrt{2\Gamma(u)(x_i)}}{\beta-u(x_i)}.
\eea
Since $\beta-u>0$ and $u\leq 1$, we have by using Theorem \ref{Theorem 1},
\bna
\f{\sqrt{2\Gamma(u)}}{\beta-u}&=&\f{\sqrt{2\Gamma(\beta-u)}}{\beta-u}\\
&\leq&\sqrt{d}\le(\f{\Delta(\beta-u)}{\beta-u}+D_\mu+\sqrt{\f{D_\mu}{d}}\ri)\\
&=&\sqrt{d}\le(\f{\lambda u}{\beta-u}+D_\mu+\sqrt{\f{D_\mu}{d}}\ri)\\
&\leq&\sqrt{d}\le(\f{1}{\beta-1}\lambda+D_\mu+\sqrt{\f{D_\mu}{d}}\ri).
\ena
This together with (\ref{uu}) implies
\be
\sum_{i=1}^{\ell}\f{|u(x_i)-u(x_{i+1})|}{\beta-u(x_i)}\leq
D{d}\le(\f{1}{\beta-1}\lambda+D_\mu+\sqrt{\f{D_\mu}{d}}\ri).\label{up}
\ee
On the other hand,
\bea\nonumber\sum_{i=1}^{\ell}\f{|u(x_i)-u(x_{i+1})|}{\beta-u(x_i)}&\geq&
\sum_{i=1}^{\ell}\log\le(1+\f{|u(x_i)-u(x_{i+1})|}{\beta-u(x_i)}\ri)\\\nonumber
&\geq&\sum_{i=1}^{\ell}\log\f{\beta-u(x_{i+1})}{\beta-u(x_i)}\\\nonumber
&=&\log\f{\beta-u(x_{\ell})}{\beta-u(x_1)}\\\label{dn}&\geq&\log\f{\beta}{\beta-1}.\eea
Combining (\ref{up}) and (\ref{dn}), we have
$$\lambda\geq (\beta-1)\le(\f{1}{Dd}\log\f{\beta}{\beta-1}-D_\mu-\sqrt{\f{D_\mu}{d}}\ri).$$
Choose $\beta$ such that $\f{1}{Dd}\log\f{\beta}{\beta-1}-D_\mu-\sqrt{\f{D_\mu}{d}}=\f{1}{Dd}$. We obtain
$$\lambda\geq\f{1}{Dd\le(\exp\le\{1+Dd\le(D_\mu+\sqrt{\f{D_\mu}{d}}\ri)\ri\}-1\ri)}.$$
This completes the proof of the theorem. $\hfill\Box$

\bigskip

 {\bf Acknowledgements.} Y. Lin is supported by the National Science Foundation of China (Grant No.11271011). Y. Yang is supported by the National Science Foundation of China (Grant No.11171347 and Grant
 No. 11471014).


\begin{thebibliography}{00}

\bibitem{BCG} M. Barlow, T. Coulhon, A. Grigor'yan, Manifolds and graphs with slow heat kernel decay, Invent. Math.
144 (2001) 609-649.

\bibitem{BHY} F. Bauer, B. Hua, S. T. Yau, Davies-Gaffney-Grior'yan lemma on graphs, Preprint.

\bibitem{BHLLY} F. Bauer, P. Horn, Y. Lin, G. Lippner, D. Mangoubi, S. T. Yau, Li-Yau inequality on graphs, J. Differential Geom.
99 (2015) 359-405.


\bibitem{Che75} S. Y. Cheng, Eigenvalue comparision theorems and its geometric applications, Math Z 143 (1975) 289-297.

\bibitem{FanChung} Fan R. K. Chung, Spectral graph theory, CBMS regional conference series in Math., No. 92, 1997, Amer. Math. Society.

\bibitem{Dav92} E. Davies, Heat kernel bounds, conservation of probability and the Feller property, J. d'Analyse Math. 58 (1992) 99-119.

\bibitem{Gaf59} M. Gaffney, The conservation property of the heat equation on Riemannian manifolds, Comm. Pure Appl. Math. 12 (1959) 1-11.

\bibitem{Gri94} A. Grigor'yan, Integral maximum principle and its applications, Proc. R. Soc. A 124 (1994) 353-362.

\bibitem{HK11} S. Haeseler, M. Keller, Generalized solutions and spectrum for Dirichlet forms on graphs, Random walks, Boundaries
and spectra, Progress in probability 64 (2011) 181-201.


\bibitem{LiJ} J. Li, Gradient estimates and Harnack inequalities for nonlinear parabolic and nonlinear elliptic
equations on Riemannian manifolds, J. Funct. Anal. 100 (1991) 233-256.

\bibitem{LiYau1} P. Li, S. T. Yau, Estimates of eigenvalues of a compact Riemannian manifold, AMS symposium on the
geometry of the Laplace operator, University of Hawaii at Manoa, 1979, 205-239.

\bibitem{LiYau2} P. Li, S. T. Yau, On the parabolic kernel of the Schr\"odinger operator, Acta Math. 156 (1986) 153-201.

\bibitem{LinYau} Y. Lin, S. T. Yau, Ricci curvature and eigenvalue estimate on locally finite graphs,
Math. Res. Lett. 17 (2010) 343-356.

\bibitem{Ma} L. Ma, Gradient estimates for a simple elliptic equation on non-compact Riemannian manifolds, J. Funct. Anal.
241 (2006) 374-382.


\bibitem{Negrin} E. Negrin, Gradient estimates and a Liouville type theorem for the Schr\"odinger operator,
J. Funct. Anal. 127 (1995) 198-203.

\bibitem{Yang-1} Y. Yang, Gradient estimates for a nonlinear parabolic equation on Riemannian manifolds,
Proc. Amer. Math. Soc. 136 (2008) 4095-4102.

\bibitem{Yang-2} Y. Yang, Gradient estimates for the equation $\Delta u+cu^{-\alpha}=0$ on Riemannian manifolds,
Acta Math. Sinica, English Series, 26 (2010) 1177-1182.




\end{thebibliography}
\end{document}